\documentclass{amsart}
\usepackage{amsmath}
\usepackage{amsfonts}

\newtheorem{theorem}{Theorem}
\theoremstyle{plain}
\newtheorem{acknowledgement}{Acknowledgement}

\newtheorem{example}{Example}

\newtheorem{lemma}{Lemma}

\input{tcilatex}

\begin{document}
\title[Inverse nodal problem]{Inverse nodal problems for Dirac-type
integro-differential operators}
\author{Baki Keskin}
\author{A. Sinan Ozkan}
\curraddr{Department of Mathematics, Faculty of Science, Cumhuriyet
University 58140 \\
Sivas, TURKEY}
\email{bkeskin@cumhuriyet.edu.tr}
\curraddr{Department of Mathematics, Faculty of Science, Cumhuriyet
University 58140 \\
Sivas, TURKEY}
\email{sozkan@cumhuriyet.edu.tr}
\subjclass[2000]{34A55, 34L05, 34K29, 34K10, 47G20}
\keywords{Dirac operator, integral-differential operators, inverse nodal
problem, uniqueness theorem.}

\begin{abstract}
The inverse nodal problem for Dirac differential operator perturbated by a
Volterra integral operator is studied. We prove that dense subset of the
nodal points determines\ the coefficients of differential and integral part\
of the operator. We also provide a uniqueness theorem and an algorithm to
reconstruct the coefficients of the problem by using the nodal points.
\end{abstract}

\maketitle

\section{\textbf{Introduction}}

We consider the boundary value problem $L$ generated by the system of Dirac
integro-differential equations:%
\begin{equation}
\ell \left[ Y(x)\right] :=BY^{\prime }(x)+\Omega
(x)Y(x)+\int\limits_{0}^{x}M(x,t)Y(t)dt=\lambda Y(x),\text{ \ }x\in (0,\pi ),
\end{equation}%
subject to the boundary conditions

\begin{eqnarray}
U(y) &:&=y_{1}(0)\sin \alpha +y_{2}(0)\cos \alpha =0\bigskip \\
V(y) &:&=y_{1}(\pi )\sin \beta +y_{2}(\pi )\cos \beta =0
\end{eqnarray}%
where $\alpha $ , $\beta $ are real constants and $\lambda $ is the spectral
parameter, $B=\left( 
\begin{array}{cc}
0 & 1 \\ 
-1 & 0%
\end{array}%
\right) ,$ \ $\Omega (x)=\left( 
\begin{array}{cc}
V(x)+m & 0 \\ 
0 & V(x)-m%
\end{array}%
\right) ,$ $M(x,t)=\left( 
\begin{array}{cc}
M_{11}(x,t) & M_{12}(x,t) \\ 
M_{21}(x,t) & M_{22}(x,t)%
\end{array}%
\right) ,$ \ $Y(x)=\left( 
\begin{array}{c}
y_{1}(x) \\ 
y_{2}(x)%
\end{array}%
\right) ,$ $\Omega (x)$ and $M\left( x,t\right) $ are real-valued functions
in the class of $W_{2}^{1}(0,\pi )$, where $m$ is a real constant.
Throughout this paper, we denote $p(x)=V(x)+m,$ $r(x)=V(x)-m$ .

In 1988, the first results of the inverse nodal Sturm--Liouville problem was
given by McLaughlin \cite{mc1} who proved that the potential of the
Sturm--Liouville problem can be determined by a given dense subset of nodal
points of the eigenfunctions. In 1989, Hald and McLaughlin consider more
general boundary conditions and give some numerical schemes for the
reconstruction of the potential from nodal points \cite{H}. Yang provided an
algorithm to solve inverse nodal Sturm--Liouville problem in 1997 \cite{yang}%
. Inverse nodal problems for Sturm--Liouville or diffusion operators have
been studied in\ the several papers (\cite{Br2}, \cite{Buterin}, \cite%
{buterin2}, \cite{ch}, \cite{Law}, \cite{Ozkan}, \cite{Yur}, \cite{Yang3}
and \cite{Yang4}). The inverse nodal problems for Dirac operators with
various boundary conditions have been solved in \cite{Guo}, \cite{Yang2} and 
\cite{Yang5}. In their works, it was shown that the zeros of the first
components of the eigenfunctions determines the coefficients of operator.

Nowadays, the studies concerning the perturbation of a differential operator
by a Volterra integral operator, namely the integro-differential operator,
are beginning to have a significant place in the literature (\cite{But}, 
\cite{But 2}, \cite{G}, \cite{Kur2} and \cite{B}). For Sturm-Liouville type
integro-differential operators, there exist some studies about inverse
problems but there is no study for Dirac type integro-differential
operators. The inverse nodal problem for Sturm-Liouville type
integro-differential operators was first studied by \cite{Kur}. In their
study, it is shown that the potential function can be determined by using
nodal points while the coefficient of the integral operator is known. In our
study, we prove that the integral operator can be partially determined as
well as the potential function and the other coefficients of the problem.

\section{\textbf{Preliminaries}}

Let $\varphi (x,\lambda )=\left( \varphi _{1}(x,\lambda ),\varphi
_{2}(x,\lambda )\right) ^{T}$ be the solution of (1) satisfying the initial
condition $\varphi (0,\lambda )=(\cos \alpha ,-\sin \alpha )^{T}$. For each
fixed $x$ and $t,$ this solution is an entire function of $\lambda .$

It is clear that $\varphi (x,\lambda )$ satisfies the following integral
equations:%
\begin{equation}
\begin{array}{l}
\varphi _{1}(x,\lambda )=\cos (\lambda x-\alpha )+\int_{0}^{x}\sin \lambda
(x-t)p(t)\varphi _{1}(t,\lambda )dt\bigskip \\ 
+\int_{0}^{x}\cos \lambda (x-t)r(t)\varphi _{2}(t,\lambda )dt\bigskip \\ 
+\int_{0}^{x}\int_{0}^{t}\sin \lambda (x-t)\left\{ M_{11}(t,\xi )\varphi
_{1}(\lambda ,\xi )+M_{12}(t,\xi )\varphi _{2}(\lambda ,\xi )\right\} d\xi
dt\bigskip \\ 
+\int_{0}^{x}\int_{0}^{t}\cos \lambda (x-t)\left\{ M_{21}(t,\xi )\varphi
_{1}(\lambda ,\xi )+M_{22}(t,\xi )\varphi _{2}(\lambda ,\xi )\right\} d\xi
dt\bigskip%
\end{array}%
\end{equation}%
\begin{equation}
\begin{array}{l}
\varphi _{2}(x,\lambda )=\sin (\lambda x-\alpha )-\int_{0}^{x}\cos \lambda
(x-t)p(t)\varphi _{1}(t,\lambda )dt\bigskip \\ 
+\int_{0}^{x}\sin \lambda (x-t)r(t)\varphi _{2}(t,\lambda )dt\bigskip \\ 
-\int_{0}^{x}\int_{0}^{t}\cos \lambda (x-t)\left\{ M_{11}(t,\xi )\varphi
_{1}(\lambda ,\xi )+M_{12}(t,\xi )\varphi _{2}(\lambda ,\xi )\right\} d\xi
dt\bigskip \\ 
+\int_{0}^{x}\int_{0}^{t}\sin \lambda (x-t)\left\{ M_{21}(t,\xi )\varphi
_{1}(\lambda ,\xi )+M_{22}(t,\xi )\varphi _{2}(\lambda ,\xi )\right\} d\xi dt%
\end{array}%
\end{equation}%
To apply the method of successive approximations to the system (4) and (5),
we denote 
\begin{equation*}
\begin{array}{l}
\varphi _{1,0}(x,\lambda )=\cos (\lambda x-\alpha ),\bigskip \\ 
\varphi _{1,n+1}(x,\lambda )=\int_{0}^{x}\sin \lambda (x-t)p(t)\varphi
_{1,n}(t,\lambda )dt\bigskip \\ 
+\int_{0}^{x}\cos \lambda (x-t)r(t)\varphi _{2,n}(t,\lambda )dt\bigskip \\ 
+\int_{0}^{x}\int_{0}^{t}\sin \lambda (x-t)\left\{ M_{11}(t,\xi )\varphi
_{1,n}(\lambda ,\xi )+M_{12}(t,\xi )\varphi _{2,n}(\lambda ,\xi )\right\}
d\xi dt\bigskip \\ 
+\int_{0}^{x}\int_{0}^{t}\cos \lambda (x-t)\left\{ M_{21}(t,\xi )\varphi
_{1,n}(\lambda ,\xi )+M_{22}(t,\xi )\varphi _{2,n}(\lambda ,\xi )\right\}
d\xi dt\bigskip%
\end{array}%
\end{equation*}%
and%
\begin{equation*}
\begin{array}{l}
\varphi _{2,0}(x,\lambda )=\sin (\lambda x-\alpha ),\bigskip \\ 
\varphi _{2,n+1}(x,\lambda )=-\int_{0}^{x}\cos \lambda (x-t)p(t)\varphi
_{1,n}(t,\lambda )dt\bigskip \\ 
+\int_{0}^{x}\sin \lambda (x-t)r(t)\varphi _{2,n}(t,\lambda )dt\bigskip \\ 
-\int_{0}^{x}\int_{0}^{t}\cos \lambda (x-t)\left\{ M_{11}(t,\xi )\varphi
_{1,n}(\lambda ,\xi )+M_{12}(t,\xi )\varphi _{2,n}(\lambda ,\xi )\right\}
d\xi dt\bigskip \\ 
+\int_{0}^{x}\int_{0}^{t}\sin \lambda (x-t)\left\{ M_{21}(t,\xi )\varphi
_{1,n}(\lambda ,\xi )+M_{22}(t,\xi )\varphi _{2,n}(\lambda ,\xi )\right\}
d\xi dt.%
\end{array}%
\end{equation*}%
Then we have%
\begin{eqnarray*}
\varphi _{1,1}(x,\lambda ) &=&\omega (x)\sin (\lambda x-\alpha )+\dfrac{%
m\sin \alpha }{\lambda }\sin \lambda x \\
&&-\dfrac{K(x)}{2\lambda }\cos (\lambda x-\alpha )-\dfrac{L(x)}{2\lambda }%
\sin (\lambda x-\alpha )+o\left( \dfrac{e^{\left\vert \tau \right\vert x}}{%
\lambda }\right) ,
\end{eqnarray*}%
\begin{eqnarray*}
\varphi _{2,1}(x,\lambda ) &=&-\omega (x)\cos (\lambda x-\alpha )-\dfrac{%
m\cos \alpha }{\lambda }\sin \lambda x \\
&&-\dfrac{K(x)}{2\lambda }\sin (\lambda x-\alpha )+\dfrac{L(x)}{2\lambda }%
\cos (\lambda x-\alpha )+o\left( \dfrac{e^{\left\vert \tau \right\vert x}}{%
\lambda }\right) ,
\end{eqnarray*}%
and for $n\geq 1$%
\begin{eqnarray*}
\varphi _{1,2n+1}(x,\lambda ) &=&\left( -1\right) ^{n}\frac{\omega ^{2n+1}(x)%
}{\left( 2n+1\right) !}\sin (\lambda x-\alpha )+\left( -1\right) ^{n}\frac{%
\omega ^{2n}(x)}{\left( 2n\right) !}\dfrac{m\sin \alpha }{\lambda }\sin
\lambda x\bigskip \\
&&+\left( -1\right) ^{n}\frac{\omega ^{2n-1}(x)}{\left( 2n-1\right) !}\dfrac{%
m^{2}x}{2\lambda }\cos (\lambda x-\alpha )-\left( -1\right) ^{n}\frac{\omega
^{2n}(x)}{\left( 2n\right) !}\dfrac{K(x)}{2\lambda }\cos (\lambda x-\alpha
)\bigskip \\
&&-\left( -1\right) ^{n}\frac{\omega ^{2n}(x)}{\left( 2n\right) !}\dfrac{L(x)%
}{2\lambda }\sin (\lambda x-\alpha )+o\left( \dfrac{e^{\left\vert \tau
\right\vert x}}{\lambda }\right) ,
\end{eqnarray*}%
\begin{eqnarray*}
\varphi _{1,2n}(x,\lambda ) &=&\left( -1\right) ^{n}\frac{\omega ^{2n}(x)}{%
\left( 2n\right) !}\cos (\lambda x-\alpha )+\left( -1\right) ^{n}\frac{%
\omega ^{2n-1}(x)}{\left( 2n-1\right) !}\dfrac{m\sin \alpha }{\lambda }\cos
\lambda x\bigskip \\
&&-\left( -1\right) ^{n}\frac{\omega ^{2n-2}(x)}{\left( 2n-2\right) !}\dfrac{%
m^{2}x}{2\lambda }\sin (\lambda x-\alpha )+\left( -1\right) ^{n}\frac{\omega
^{2n-1}(x)}{\left( 2n-1\right) !}\dfrac{K(x)}{2\lambda }\sin (\lambda
x-\alpha )\bigskip \\
&&-\left( -1\right) ^{n}\frac{\omega ^{2n-1}(x)}{\left( 2n-1\right) !}\dfrac{%
L(x)}{2\lambda }\cos (\lambda x-\alpha )+o\left( \dfrac{e^{\left\vert \tau
\right\vert x}}{\lambda }\right) ,
\end{eqnarray*}%
\begin{eqnarray*}
\varphi _{2,2n+1}(x,\lambda ) &=&-\left( -1\right) ^{n}\frac{\omega
^{2n+1}(x)}{\left( 2n+1\right) !}\cos (\lambda x-\alpha )-\left( -1\right)
^{n}\frac{\omega ^{2n}(x)}{\left( 2n\right) !}\dfrac{m\cos \alpha }{\lambda }%
\sin \lambda x\bigskip \\
&&+\left( -1\right) ^{n}\frac{\omega ^{2n-1}(x)}{\left( 2n-1\right) !}\dfrac{%
m^{2}x}{2\lambda }\sin (\lambda x-\alpha )-\left( -1\right) ^{n}\frac{\omega
^{2n}(x)}{\left( 2n\right) !}\dfrac{K(x)}{2\lambda }\sin (\lambda x-\alpha
)\bigskip \\
&&+\left( -1\right) ^{n}\frac{\omega ^{2n}(x)}{\left( 2n\right) !}\dfrac{L(x)%
}{2\lambda }\cos (\lambda x-\alpha )+o\left( \dfrac{e^{\left\vert \tau
\right\vert x}}{\lambda }\right) ,
\end{eqnarray*}%
\begin{eqnarray*}
\varphi _{2,2n}(x,\lambda ) &=&\left( -1\right) ^{n}\frac{\omega ^{2n}(x)}{%
\left( 2n\right) !}\sin (\lambda x-\alpha )-\left( -1\right) ^{n}\frac{%
\omega ^{2n-1}(x)}{\left( 2n-1\right) !}\dfrac{m\cos \alpha }{\lambda }\cos
\lambda x\bigskip \\
&&+\left( -1\right) ^{n}\frac{\omega ^{2n-2}(x)}{\left( 2n-2\right) !}\dfrac{%
m^{2}x}{2\lambda }\cos (\lambda x-\alpha )-\left( -1\right) ^{n}\frac{\omega
^{2n-1}(x)}{\left( 2n-1\right) !}\dfrac{K(x)}{2\lambda }\cos (\lambda
x-\alpha )\bigskip \\
&&-\left( -1\right) ^{n}\frac{\omega ^{2n-1}(x)}{\left( 2n-1\right) !}\dfrac{%
L(x)}{2\lambda }\sin (\lambda x-\alpha )+o\left( \dfrac{e^{\left\vert \tau
\right\vert x}}{\lambda }\right) ,
\end{eqnarray*}%
where, $\omega (x)=\dfrac{1}{2}\int_{0}^{x}(p(t)+r(t))dt=\int_{0}^{x}V(t)dt,$
$K(x)=\int_{0}^{x}(M_{11}(t,t)+M_{22}(t,t))dt,$ $L(x)=%
\int_{0}^{x}(M_{12}(t,t)-M_{21}(t,t))dt$ and $\tau =\func{Im}\lambda .$

Thus, the functions $\varphi _{1}(x,\lambda )$ and $\varphi _{2}(x,\lambda )$
have the following asymptotic formulae:%
\begin{eqnarray}
\varphi _{1}(x,\lambda ) &=&\cos (\lambda x-\omega (x)-\alpha )+\dfrac{m\sin
\alpha }{\lambda }\sin (\lambda x-\omega (x))\bigskip \\
&&+\dfrac{m^{2}x}{2\lambda }\sin (\lambda x-\omega (x)-\alpha )-\dfrac{K(x)}{%
2\lambda }\cos (\lambda x-\omega (x)-\alpha )\bigskip  \notag \\
&&-\dfrac{L(x)}{2\lambda }\sin (\lambda x-\omega (x)-\alpha )+o\left( \dfrac{%
e^{\left\vert \tau \right\vert x}}{\lambda }\right) ,  \notag
\end{eqnarray}%
\begin{eqnarray}
\varphi _{2}(x,\lambda ) &=&\sin (\lambda x-\omega (x)-\alpha )-\dfrac{m\cos
\alpha }{\lambda }\sin (\lambda x-\omega (x))\bigskip \\
&&-\dfrac{m^{2}x}{2\lambda }\cos (\lambda x-\omega (x)-\alpha )-\dfrac{K(x)}{%
2\lambda }\sin (\lambda x-\omega (x)-\alpha )\bigskip  \notag \\
&&+\dfrac{L(x)}{2\lambda }\cos (\lambda x-\omega (x)-\alpha )+o\left( \dfrac{%
e^{\left\vert \tau \right\vert x}}{\lambda }\right)  \notag
\end{eqnarray}%
for sufficiently large $\left\vert \lambda \right\vert ,$ uniformly in $x.$%
\bigskip

The characteristic function $\Delta (\lambda )$ of the problem (1)-(3) is
defined by the relation 
\begin{equation}
\Delta (\lambda )=\varphi _{1}(\pi ,\lambda )\sin \beta +\varphi _{2}(\pi
,\lambda )\cos \beta ,
\end{equation}%
It is obvious that $\Delta (\lambda )$ is an entire function and its zeros,
namely $\left\{ \lambda _{n}\right\} _{n\geq 0}$ ,\ coincide with the
eigenvalues of the problem (1)-(3). Using the asymptotic formulae (6) and
(7), one can easily obtain%
\begin{eqnarray}
\Delta (\lambda ) &=&\sin (\lambda \pi -\omega (\pi )+\beta -\alpha )-\dfrac{%
m^{2}\pi }{2\lambda }\cos (\lambda \pi -\omega (\pi )+\beta -\alpha ) \\
&&-\dfrac{K(\pi )}{2\lambda }\sin (\lambda \pi -\omega (\pi )+\beta -\alpha
)+\dfrac{L(\pi )}{2\lambda }\cos (\lambda \pi -\omega (\pi )+\beta -\alpha )
\notag \\
&&-\dfrac{m}{\lambda }\sin (\lambda \pi -\omega (\pi ))\cos (\beta +\alpha
)+o\left( \dfrac{e^{\left\vert \tau \right\vert \pi }}{\lambda }\right) 
\notag
\end{eqnarray}%
\bigskip for sufficiently large $\left\vert \lambda \right\vert .$ Since the
eigenvalues of the problem (1)-(3) are the roots of $\Delta (\lambda _{n})=0$%
, we can write the following equation for them:$\bigskip $\newline
$\left( 1-\dfrac{K(\pi )}{2\lambda _{n}}-\dfrac{m}{\lambda _{n}}\cos (\beta
+\alpha )\cos (\alpha -\beta )\right) \tan (\lambda _{n}\pi -\omega (\pi
)+\beta -\alpha )\bigskip $

$\bigskip =\dfrac{m^{2}\pi }{2\lambda _{n}}-\dfrac{L(\pi )}{2\lambda _{n}}+%
\dfrac{m}{\lambda _{n}}\cos (\beta +\alpha )\sin (\alpha -\beta )+o\left( 
\dfrac{e^{\left\vert \tau \right\vert \pi }}{\lambda _{n}}\right) $\newline
which implies that$\bigskip $

\bigskip $\tan (\lambda _{n}\pi -\omega (\pi )+\beta -\alpha )=\left( \dfrac{%
m^{2}\pi }{2\lambda _{n}}-\dfrac{L(\pi )}{2\lambda _{n}}+\dfrac{m}{\lambda
_{n}}\cos (\beta +\alpha )\sin (\alpha -\beta )\right) \times $

$\ \ \ \ \ \ \ \ \ \ \ \ \ \ \ \ \ \ \ \ \ \ \ \ \ \ \ \ \ \ \ \ \ \ \ \ \ \
\ \ \ \times \left( 1+\dfrac{K(\pi )}{2\lambda _{n}}+\dfrac{m}{\lambda _{n}}%
\cos (\beta +\alpha )\cos (\alpha -\beta )+o\left( \dfrac{1}{\lambda _{n}}%
\right) \right) \bigskip $\newline
\bigskip for sufficiently large $n.$ We obtain from the last equation,

\begin{eqnarray}
\lambda _{n} &=&n+\dfrac{1}{\pi }\int_{0}^{\pi }V(t)dt+\dfrac{\alpha -\beta 
}{\pi }\bigskip \\
&&+\dfrac{1}{2n\pi }\left( m^{2}\pi -L(\pi )+2m\cos (\beta +\alpha )\sin
(\alpha -\beta )\right) +o\left( \dfrac{1}{n}\right)  \notag
\end{eqnarray}%
for $n=0,\pm 1,\pm 2,...$

\section{\textbf{Main Results}}

\begin{lemma}
For sufficiently large $n$, the first component $\varphi _{1}(x,\lambda
_{n}) $ of the eigenfunction $\varphi (x,\lambda _{n})$ has exactly $n$
nodes $\left\{ x_{n}^{j}:j=\overline{0,n}\right\} $ in the interval $\left(
0,\pi \right) $:\newline
$0<x_{n}^{0}<x_{n}^{1}<...<x_{n}^{n}<\pi $. The numbers $\left\{
x_{n}^{j}\right\} $ satisfy the following asymptotic formula:
\end{lemma}

\begin{equation}
\begin{array}{l}
x_{n}^{j}=\dfrac{\left( j+1/2\right) \pi }{n}+\dfrac{\omega
(x_{n}^{j})+\alpha }{n}\bigskip \\ 
\text{ \ \ }-\dfrac{\left( j+1/2\right) \pi }{n}\left( \dfrac{\omega (\pi
)+\alpha -\beta }{n\pi }\right) -\dfrac{\omega (\pi )+\alpha -\beta }{%
n^{2}\pi }\left( \omega (x_{n}^{j})+\alpha \right) \\ 
\text{ \ \ }+\dfrac{1}{2n^{2}}\left( m^{2}x_{n}^{j}-L(x_{n}^{j})+m\sin
2\alpha \right) +o\left( \frac{1}{n^{2}}\right) .%
\end{array}%
\end{equation}

\begin{proof}
From (6), the following asymptotic formula can be written for sufficiently
large $n.\bigskip $\newline
$\varphi _{1}(x,\lambda _{n})=\cos (\lambda _{n}x-\omega (x)-\alpha )+\dfrac{%
m\sin 2\alpha }{2\lambda _{n}}\sin (\lambda _{n}x-\omega (x)-\alpha )$

$+\dfrac{m\sin ^{2}\alpha }{\lambda _{n}}\cos (\lambda _{n}x-\omega
(x)-\alpha )+\dfrac{m^{2}x}{2\lambda _{n}}\sin (\lambda _{n}x-\omega
(x)-\alpha )$

$-\dfrac{K(x)}{2\lambda _{n}}\cos (\lambda _{n}x-\omega (x)-\alpha )-\dfrac{%
L(x)}{2\lambda _{n}}\sin (\lambda _{n}x-\omega (x)-\alpha )+o\left( \dfrac{%
e^{\left\vert \tau _{n}\right\vert x}}{\lambda _{n}}\right) \bigskip $%
\newline
From $\varphi _{1}(x_{n}^{j},\lambda _{n})=0,$ we get$\bigskip $\newline
$\cos \left( \lambda _{n}x_{n}^{j}-\omega (x_{n}^{j})-\alpha \right) =\dfrac{%
m\sin 2\alpha }{2\lambda _{n}}\sin (\lambda _{n}x_{n}^{j}-\omega
(x_{n}^{j})-\alpha )\medskip $

$+\dfrac{m\sin ^{2}\alpha }{\lambda _{n}}\cos (\lambda _{n}x_{n}^{j}-\omega
(x_{n}^{j})-\alpha )\medskip $

$+\dfrac{m^{2}x_{n}^{j}}{2\lambda _{n}}\sin (\lambda _{n}x_{n}^{j}-\omega
(x_{n}^{j})-\alpha )\medskip $

$\medskip -\dfrac{K(x_{n}^{j})}{2\lambda _{n}}\cos (\lambda
_{n}x_{n}^{j}-\omega (x_{n}^{j})-\alpha )$

$-\dfrac{L(x_{n}^{j})}{2\lambda _{n}}\sin (\lambda _{n}x_{n}^{j}-\omega
(x_{n}^{j})-\alpha )+o\left( \dfrac{e^{\left\vert \tau _{n}\right\vert x}}{%
\lambda _{n}}\right) .\medskip $\newline
which is equivalent to

$\tan \left( \lambda _{n}x_{n}^{j}-\omega (x_{n}^{j})-\alpha -\dfrac{\pi }{2}%
\right) =\left( 1-\dfrac{K(x)}{2\lambda _{n}}+\dfrac{m\sin ^{2}\alpha }{%
\lambda _{n}}\right) ^{-1}\times $

$\ \ \ \ \ \ \ \ \ \ \ \ \ \ \ \ \ \ \ \ \ \ \ \ \ \ \ \ \ \ \ \ \ \ \ \ \ \
\ \ \times \left( \dfrac{m^{2}x_{n}^{j}}{2\lambda _{n}}-\dfrac{L(x_{n}^{j})}{%
2\lambda _{n}}+\dfrac{m\sin 2\alpha }{\lambda _{n}}+o\left( \dfrac{1}{%
\lambda _{n}}\right) \right) .$\newline
Taking into account Taylor's expansion formula for the arctangent, we get$%
\medskip $ 
\begin{equation*}
\lambda _{n}x_{n}^{j}-\omega (x_{n}^{j})-\alpha -\dfrac{\pi }{2}=j\pi +%
\dfrac{1}{2\lambda _{n}}\left( m^{2}x_{n}^{j}-L(x_{n}^{j})+m\sin 2\alpha
\right) +o\left( \dfrac{1}{\lambda _{n}}\right) .\medskip
\end{equation*}%
It follows from the last equality

\begin{equation*}
x_{n}^{j}=\dfrac{\left( j+\frac{1}{2}\right) \pi +\omega (x_{n}^{j})+\alpha 
}{\lambda _{n}}+\dfrac{1}{2\lambda _{n}^{2}}\left(
m^{2}x_{n}^{j}-L(x_{n}^{j})+m\sin 2\alpha \right) +o\left( \dfrac{1}{\lambda
_{n}^{2}}\right) .
\end{equation*}%
The relation (11) is proven by using the asymptotic formula%
\begin{equation*}
\lambda _{n}^{-1}=\dfrac{1}{n}\left\{ 1-\dfrac{\omega (\pi )+\alpha -\beta }{%
n\pi }-\dfrac{\left( m^{2}\pi -L(\pi )+2m\cos (\beta +\alpha )\sin (\alpha
-\beta )\right) }{2n^{2}\pi }+o\left( \dfrac{1}{n^{2}}\right) \right\}
\end{equation*}%
From (11), we have the following asymptotic expression for nodal lenghts:

\begin{equation*}
l_{n}^{j}:=x_{n}^{j+1}-x_{n}^{j}=\dfrac{\pi }{n}+o\left( \dfrac{1}{n}\right)
.
\end{equation*}%
One can easily see that $\varphi _{1}(\dfrac{k\pi }{n},\lambda _{n})$ and $%
\varphi _{1}(\dfrac{\left( k+1\right) \pi }{n},\lambda _{n})$ have different
signs for each fixed $k$ and for sufficiently large $n.$ Thus, the function $%
\varphi _{1}(x,\lambda _{n})$ has exactly $n$ nodes in $\left( 0,\pi \right)
.$
\end{proof}

Let $X$ be the set of nodal points and $\omega (\pi )=0.$ For each fixed $%
x\in \left( 0,\pi \right) $ $\ $we can choose a sequence $\left(
x_{n}^{j}\right) \subset X$ so that $x_{n}^{j}$ converges to $x.$ Then the
following limits are exist and finite:%
\begin{equation}
\underset{n\rightarrow \infty }{\lim }n\left( x_{n}^{j(n)}-\frac{\left( j+%
\frac{1}{2}\right) \pi }{n}\right) =f(x),
\end{equation}%
where%
\begin{equation*}
f(x)=\omega (x)+\alpha -\dfrac{x}{\pi }\left( \alpha -\beta \right)
\end{equation*}%
and%
\begin{equation}
\underset{n\rightarrow \infty }{\lim }2n^{2}\left( x_{n}^{j}-\frac{\left( j+%
\frac{1}{2}\right) \pi -\omega (x_{n}^{j})}{n}+\dfrac{\left( j+1/2\right)
\pi }{n}\left( \dfrac{\alpha -\beta }{n\pi }\right) \right) =g(x),
\end{equation}%
where%
\begin{equation*}
g(x)=-L(x)+2\dfrac{\beta -\alpha }{\pi }\left( \omega (x)+\alpha \right)
+m^{2}x+m\sin 2\alpha
\end{equation*}%
Therefore, proof of the following theorem is clear.

\begin{theorem}
The given dense subset of nodal points $X$ uniquely determines the potential 
$V(x)$, the function $L^{\prime }(x)=M_{12}(x,x)-M_{21}(x,x)$ a.e. on $%
\left( 0,\pi \right) ,$ and the coefficients $\alpha $ and $\beta $ of the
boundary conditions. Moreover, $V(x),$ $L^{\prime }(x),$ $\alpha $ and $%
\beta $ can be reconstructed by the following formulae:

\textbf{Step-1:} For each fixed $x\in (0,\pi ),$ choose a sequence $\left(
x_{n}^{j(n)}\right) \subset X$ such that $\underset{n\rightarrow \infty }{%
\lim }x_{n}^{j(n)}=x;$

\textbf{Step-2: }Find the function $f(x)$ from (12) and calculate 
\begin{eqnarray*}
V(x) &=&\medskip f^{\prime }(x) \\
\alpha &=&f(0) \\
\beta &=&f(\pi )
\end{eqnarray*}

\textbf{Step-3: } If $\sin 2\alpha \neq 0,$ find the function $g(x)$ from
(13) and calculate%
\begin{eqnarray*}
\medskip m &=&\dfrac{g(0)+2\alpha \left( \alpha -\beta \right) }{\sin
2\alpha } \\
L^{\prime }(x) &=&\medskip -g^{\prime }(x)+2\dfrac{\beta -\alpha }{\pi }%
V(x)+m^{2}
\end{eqnarray*}%
otherwise assume $m$ is known.
\end{theorem}

\begin{example}
Let $\left\{ x_{n}^{j}\right\} \subset X$ be the dense subset of nodal
points in $(0,\pi )$ given by the following asimptotics:\newline
$x_{n}^{j}=\dfrac{\left( j+1/2\right) \pi }{n}+\dfrac{\frac{\pi }{4}+\sin 
\dfrac{\left( j+1/2\right) \pi }{n}}{n}$

$+\dfrac{1}{2n^{2}}\left( \dfrac{\left( j+1/2\right) \pi }{n}+\sin \dfrac{%
\left( j+1/2\right) \pi }{n}+1\right) +o\left( \frac{1}{n^{2}}\right)
.\bigskip $\newline
It can \ be calculated from (12) and (13) that,%
\begin{eqnarray*}
f(x) &=&\frac{\pi }{4}+\sin x \\
g(x) &=&1+x+\sin x
\end{eqnarray*}%
Therefore, it is obtained \ by using the algorithm in Therem 1,%
\begin{eqnarray*}
&&V(x)=f^{\prime }(x)=\cos x,\medskip  \\
&&\left. \alpha =f(0)=\frac{\pi }{4}=f(\pi )=\beta ,\medskip \right.  \\
&&\left. M_{12}(x,x)-M_{21}(x,x)=L^{\prime }(x)=-\cos x\right.  \\
&&\left. m=1.\right. 
\end{eqnarray*}

\begin{acknowledgement}
The similar results can be obtained for the equation (1) with parameter
dependent boundary conditions. 
\end{acknowledgement}
\end{example}

\end{document}